\documentclass[11pt,reqno]{amsart}
\usepackage{amssymb,amscd,amsbsy}
\usepackage{amssymb,amscd,amsbsy,mathrsfs}
\newcommand{\lb}{\linebreak}

\renewcommand{\a}{\alpha}

\newcommand{\z}{\zeta}

\newcommand{\f}{\varphi}

\newcommand{\D}{\Delta}
\renewcommand{\L}{\Lambda}

\newcommand{\h}{{\mathscr H}}

\newcommand{\X}{{\mathscr X}}
\newcommand{\Y}{{\mathscr Y}}

\newcommand{\C}{{\Bbb C}}
\newcommand{\T}{{\Bbb T}}

\newcommand{\dd}{{\Bbb D}}
\newcommand{\R}{{\Bbb R}}

\newcommand{\bs}{\boldsymbol}

\newcommand{\m}{{\boldsymbol m}}
\newcommand{\bS}{{\boldsymbol S}}

\newcommand{\rf}[1]{(\ref{#1})}

\newcommand{\df}{\stackrel{\mathrm{def}}{=}}

\newcommand{\clos}{\operatorname{clos}}

\newcommand{\const}{\operatorname{const}}

\newcommand{\eeq}{\end{equation}}
\newcommand{\beq}{\begin{equation}}
\newcommand{\bay}{\begin{eqnarray}}
\newcommand{\ba}{\begin{align*}}
\newcommand{\ea}{\end{align*}}
\newcommand{\ey}{\end{eqnarray}}
\newcommand{\bey}{\begin{eqnarray*}}
\newcommand{\eey}{\end{eqnarray*}}

\newcommand{\be}{\infty}

\newcommand{\bl}{\blacksquare}

\newcommand{\Pf}{{\bf Proof. }}
\newcommand{\im}{\operatorname{Im}}

\newcommand{\ov}{\overline}

\newtheorem{thm}{\hspace{\parindent}Theorem}[section]

\newtheorem{lem}[thm]{\hspace{\parindent}Lemma}

\theoremstyle{remark}

\newtheorem*{rem*}{Remark}

\newcommand\mS{\mathcal{S}}
\newcommand\mT{\mathcal{T}}

\newcommand\fM{\frak M}

\newcommand{\rd}{{\rm d}}

\newcommand\mB{\mathcal{B}}

\newcommand\ri{{\rm i}}

\newcommand{\CAr}{{\rm C}_{\rm A}(\C_+)}

\newcommand{\CAb}{{\rm C}_{{\rm A},\be}}

\begin{document}

\author{A.B. Aleksandrov and V.V. Peller}

\title{Analytic Schur multipliers}
\thanks{The research on \S\:2 is supported by 
Russian Science Foundation [grant number 23-11-00153].
The research on on \S\:3 is supported by a grant of the Government of the Russian Federation for the state support of scientific research, carried out under the supervision of leading scientists, agreement  075-15-2025-013.}

\begin{abstract}
We study in this paper analytic Schur multipliers on $\C_+^2$ and $\dd^2$, i.e. Schur multipliers on $\R^2$ and $\T^2$ that are boundary-value functions of functions analytic in $\C_+^2$ and $\dd^2$. Such Schur multipliers are important when studying properties of functions of maximal dissipative operators and contractions under perturbation. We show that if the boundary-value function of a Schur multiplier has certain regularity properties, then it can be represented as an element of the Haagerup tensor product of spaces of analytic functions with similar regularity properties. 
\end{abstract}

\maketitle

\numberwithin{equation}{section}

\

\setcounter{section}{0}
\section{\bf Introduction}
\label{In}
\setcounter{equation}{0}

\

This paper can be considered as a continuation of our previous work on Schur multipliers, see \cite{PeHo}, \cite{APol} and \cite{APHtp}. The notion of Schur multipliers with respect to spectral measures was introduced in
\cite{PeHo} in connection with double operator integrals; it is a natural generalization of the notion of matrix Schur multipliers, see \cite{Sch}. 

Recall that double operator integrals are expressions of the form
\bay
\label{dvopi}
\iint_{\X\times\Y}\Phi(x,y)\,\rd E_1(x)Q\,\rd E_2(y).
\ey
Here $E_1$ and $E_2$ are spectral measures on Hilbert space, $\Phi$ is a bounded measurable function and 
$Q$ is a bounded linear operator on Hilbert space. 

Double operator integrals appeared in the paper \cite{DK} by Yu.L. Daletskii and S.G. Krein. Later M.Sh. Birman  and M.S. Solomyak in \cite{BSh1}  created a beautiful theory of double operator integrals. They started with the case when $Q$ is a Hilbert--Schmidt operator and defined the double operator integral in \rf{dvopi} as
$$
\left(\int_{\X\times\Y}\Phi\,\rd{\rm E}\right)Q,
$$
where ${\rm E}$ is the spectral measure on the Hilbert Schmidt class $\bS_2$, which is defined on the measurable rectangles by
$$
{\rm E}(\L\times\D)Q=E_1(\L)QE_2(\D),\quad Q\in\bS_2.
$$
Clearly,
$$
\left\|\iint\Phi(x,y)\,\rd E_1(x)Q\,\rd E_2(y)\right\|_{\bS_2}\le\|\Phi\|_{L^\be}\|Q\|_{\bS_2}.
$$

 \medskip
 
 {\bf Definition 1.}  A bounded measurable function $\Phi$ on $\X\times\Y$ is called a {\it Schur multiplier (of $\bS_1$) with respect to the spectral measures $E_1$ and $E_2$} if
 $$
 Q\in\bS_1\quad\Longrightarrow\quad\iint\Phi\,\rd E_1Q\,\rd E_2\in\bS_1.
 $$
 We denote the class of Schur multipliers of $\bS_1$ with respect to spectral measures $E_1$ and $E_2$ by $\fM_{E_1,E_2}$.
 
 \medskip
 
 If $\Phi\in\fM_{E_2,E_1}$, we can define by duality the double operator integral in \rf{dvopi} for arbitrary bounded linear operator $Q$ and the transformer
 $$
 Q\mapsto\iint\Phi\,\rd E_1Q\,\rd E_2
 $$
 becomes a bounded linear operator in the operator norm.
 
  \medskip
 
 We also use the notation $\fM(\R^2)$ for the class of {\it Borel Schur multipliers}, i.e. the class of 
 functions on $\R^2$ that are Schur multipliers with respect to arbitrary Borel spectral measure on $\R$.
 
 We refer to \cite{APol} and \cite{APHtp} for more detailed information about double operator integrals.
 
 Let us proceed now to the definition of the Haagerup tensor product of subspaces of $L^\be$ spaces.
 
  \medskip
  
 {\bf Definition 2.}
 Suppose that $X$ and $Y$ are subspaces of $L^\be$ spaces of functions on sets $\X$ and $\Y$. The {\it Haagerup tensor product} $X\otimes_{\rm h}Y$ of $X$ and $Y$ is, by definition, the class of functions $\Phi$
on $\X\times\Y$ of the form
\bay
\label{PhiHaa}
\Phi(x,y)=\sum_{n\ge0}\f_n(x)\psi_n(y),\quad\mbox{for}\quad\f_n\in X\quad\mbox{and}\quad\psi_n\in Y
\ey
such that 
\bay
\label{uslHaa}
\sum_{n\ge0}|\f_n|^2\in L^\be\quad\mbox{and}\quad\sum_{n\ge0}|\psi_n|^2\in L^\be.
\ey
The norm of $\Phi$ in $X\otimes_{\rm h}Y$, is the infimum of the expression
$$
\left\|\left(\sum_{n\ge0}|\f_n|^2\right)^{1/2}\right\|_{L^\be}\cdot\left\|\left(\sum_{n\ge0}|\psi_n|^2\right)^{1/2}\right\|_{L^\be}
$$
over all representations of $\Phi$ in the form of \rf{PhiHaa}.

\medskip

In the same way we can define the Haagerup tensor products of two subspaces of spaces of bounded continuous functions on locally compact subspaces. 

There are various descriptions of the space $\fM_{E_1,E_2}$ of Schur multipliers, see \cite{PeHo}, \cite{APol} and \cite{APHtp}. We mention the following one given in \cite{PeHo}: $\Phi\in\fM_{E_1,E_2}$ if and only 
$\Phi\in L^\be_{E_1}\otimes_{\rm h}L^\be_{E_2}$. Note that the implication
$$
\Phi\in L^\be_{E_1}\otimes_{\rm h}L^\be_{E_2}\quad\Longrightarrow\quad\Phi\in\fM_{E_1,E_2}
$$
was obtained earlier in \cite{BSh2}.

Note also that in \cite{APHtp} it was shown that
$$
\|\Phi\|_{\fM_{E_1,E_2}}=\|\Phi\|_{L^\be_{E_1}\otimes_{\rm h}L^\be_{E_2}},\quad\Phi\in\fM_{E_1,E_2}.
$$

It would be natural to expect that if a Schur multiplier $\Phi$ in $\fM(\R^2)$ has certain regularity properties, then on can select the functions $\f_n$ and $\psi_n$ in \rf{PhiHaa} that possess similar regularity properties. 

Indeed, it can easily be deduced from Theorem 2.2.3 in \cite{APol} that if $\Phi\in\fM(\R^2)$ and for every $x\in\R$, the functions
$\Phi(\cdot,x)$ and $\Phi(x,\cdot)$ belong to the space ${\rm C}_{\rm b}(\R)$ of bounded continuous functions on $\R$, then $\Phi$ admits a representation of the form \rf{PhiHaa} with $\f_n$ and $\psi_n$ in ${\rm C}_{\rm b}(\R)$ satisfying \rf{uslHaa}.

Also, it can be deduced from Theorem 2.2.3 in \cite{APol} that if $\Phi\in\fM(\R^2)$ and for every $x\in\R$, the functions $\Phi(\cdot,x)$ and $\Phi(x,\cdot)$ belong to the space ${\rm C}(\widehat\R)$ of continuous functions $f$ on $\R$ such that the limit $\lim_{|t|\to\be} f(t)$ exists, then $\Phi$ admits a representation of the form \rf{PhiHaa} with $\f_n$ and $\psi_n$ in ${\rm C}(\widehat\R)$ satisfying \rf{uslHaa}.

In this paper we are going to obtain analogs of such results for analytic Schur multipliers. 
To begin with, we introduce the space $H^\be_{\rm r}(\C_+^2)$ of functions $\Phi$ on $\clos\C_+^2$
such that $\Phi|\C_+^2\in H^\be(\C_+^2)$ and 
$$
\lim_{t\to 0^+}\Phi(z+t\ri,w)=\lim_{t\to 0^+}\Phi(z,w+t\ri)=\Phi(z,w)
$$
for arbitrary $z,w\in\clos\C_+$.

%
%
%
%
%
%
%

\medskip

{\bf Definition 3.}
We denote by $\fM_{\rm A}(\C_+^2)$ the class of {\it analytic Schur multipliers} on $\clos\C_+^2$, i.e. the class of
functions  $f\in H^\be_{\rm r}(\C_+^2)$ such that $f|\R^2$ belongs to the space $\fM(\R^2)$ of Borel Schur multipliers.

\medskip

Throughout the paper we use the same notation for a bounded analytic function in $\C_+^2$ (or in 
$\C_+^2$) and its boundary-value function on $\R^2$ (or $\R$). Here $\C_+\df\{\z\in\C:~\im\z>0\}$.

In \S\:\ref{ogrnepanal} we characterize the analytic Schur multipliers of class 
$$
{\rm C}^{\rm S}_{{\rm A},\be}(\C_+^2)\df\big\{f\in H^\be(\C_+^2):\:~f(z,\cdot),\:~ f(\cdot,z)\in{\rm C}_{{\rm A},\be}\quad\mbox{for all}\quad z\in\clos\C_+\big\},
$$
where
$$
{\rm C}_{{\rm A},\be}\df\big\{f\in H^\be(\C_+):~f\quad\mbox{is continuous on}\quad\R\big\}.
$$

Another result in \S\:\ref{ogrnepanal} is a characterization of the analytic Schur multipliers of class 
$$
{\rm C}^{\rm S}_{{\rm A}}(\C_+^2)=\big\{f\in H^\be(\C_+^2):~f(z,\cdot), ~\:f(\cdot,z)\in{\rm C}_{{\rm A}}(\C_+)\quad\mbox{for all}\quad z\in\clos\C_+\big\},
$$
where
$$
\CAr\df\Big\{f\in\CAb:~\mbox{the limit}~\lim_{|\z|\to\be}f(\z)~\:\mbox{exists}\Big\}.
$$

Note that the restrictions of functions in ${\rm C}^{\rm S}_{{\rm A},\be}(\C_+^2)$ and in 
${\rm C}^{\rm S}_{{\rm A}}(\C_+^2)$ to $\R^2$ do not have to be continuous functions on $\R^2$.
On the other hand, the restrictions of functions in ${\rm C}^{\rm S}_{{\rm A},\be}(\C_+^2)$ and in 
${\rm C}^{\rm S}_{{\rm A}}(\C_+^2)$ to $\C_+^2$ must be holomorphic in $\C_+^2$.

In \S\:\ref{vsyudu}
we obtain a representation of the space $\fM_{\rm A}(\C_+^2)$ as the Haagerup tensor product \lb
$H^\be_{\rm r}(\C_+)\otimes_{\rm h}H^\be_{\rm r}(\C_+)$, where
\bay
\label{HrbeC}
H^\be_{\rm r}(\C_+)\df \big\{f\in H^\be(\C_+):~ \lim_{t\to0^+}
f(x+t\ri)~\:\mbox{exists for all}~\: x\in\R\big\}.
\ey

Finally, in \S\:\ref{bid} we obtain similar results on analytic Schur multipliers on the bidisc.

\

\section{\bf Analytic Schur multipliers in the spaces $\bs{{\rm C}^{\rm S}_{{\rm A},\be}(\C_+^2)}$ and 
$\bs{{\rm C}^{\rm S}_{{\rm A}}(\C_+^2)}$}
\label{ogrnepanal}
\setcounter{equation}{0}

\

In this section we characterize analytic Schur multipliers in 
${\rm C}^{\rm S}_{{\rm A},\be}(\C_+^2)$ and ${\rm C}^{\rm S}_{{\rm A}}(\C_+^2)$. The results we are going to obtain are important in the study of functions of perturbed dissipative operators, see \cite{APar}.

To begin with we recall the definition of matrix (discrete) Schur multipliers.

\medskip

{\bf Definition 4.}
Let $\mS$ and $\mT$ be arbitrary nonempty sets. If $\{a(s,t)\}_{(s,t)\in \mS\times \mT}$ is a matrix that induces a bounded linear operator from $\ell^2(\mT)$ to $\ell^2(\mS)$, then by
\lb$\|\{a(s,t)\}_{(s,t)\in\mS\times \mT}\|$, we mean the operator norm of the corresponding operator.
We use the notation $\mB(\mS\times\mT)$ for the class of matrices $\{a(s,t)\}_{(s,t)\in\mS\times\mT}$ that
induce bounded operators from $\ell^2(\mT)$ to $\ell^2(\mS)$.

A matrix $\Phi=\{\Phi(s,t)\}_{(s,t)\in\mS\times\mT}$ is called a {\it matrix Schur multiplier} (in other words, a {\it discrete
Schur multiplier}) on  $\mS\times\mT$ if
$$
A=\{a(s,t)\}_{(s,t)\in\mS\times\mT}\in\mB(\mS\times\mT)\quad\Longrightarrow\quad \Phi\star A\df \{\Phi(s,t)a(s,t)\}_{(s,t)\in\mS\times\mT}\in\mB(\mS\times\mT).
$$ 
We denote by $\fM_{\rm d}(\mS\times\mT)$
the space of matrix Schur multipliers on  $\mS\times\mT$. By the norm $\|A\|_{\fM_{\rm d}(\mS\times\mT)}$
of $A$ in $\fM_{\rm d}(\mS\times\mT)$ we mean the norm of the transformer 
$$
\Phi\mapsto \Phi\star A
$$
on the space $\mB(\mS\times\mT)$.

\medskip

It is easy to see that the space of Borel Schur multipliers $\fM(\R^2)$ is contained in $\fM_{\rm d}(\R^2)$.

\begin{thm} 
\label{beskon}
$$
\fM_{\rm A}(\C_+^2)\cap{\rm C}_{{\rm A},\be}^{\rm S}(\C_+^2)={\rm C}_{{\rm A},\be}(\C_+)\otimes_{\rm h}{\rm C}_{{\rm A},\be}(\C_+)
$$
and 
$$
\|\Phi\|_{{\rm C}_{{\rm A},\be}\otimes_{\rm h}{\rm C}_{{\rm A},\be}}
=\|\Phi\|_{\fM(\R^2)}\quad
\mbox{for every}\quad \Phi\in{\rm C}_{{\rm A},\be}(\C_+)\otimes_{\rm h}{\rm C}_{{\rm A},\be}(\C_+).
$$
\end{thm}

\Pf First we prove that
$$
{\rm C}_{{\rm A},\be}(\C_+)\otimes_{\rm h}{\rm C}_{{\rm A},\be}(\C_+)\subset\fM_{\rm A}(\C_+^2)\cap{\rm C}_{{\rm A},\be}^{\rm S}(\C_+^2)
$$
and $\|\Phi\|_{\fM(\R^2)}\le\|\Phi\|_{{\rm C}_{{\rm A},\be}\otimes_{\rm h}{\rm C}_{{\rm A},\be}}$ for all
$\Phi\in{\rm C}_{{\rm A},\be}(\C_+)\otimes_{\rm h}{\rm C}_{{\rm A},\be}(\C_+)$.

It is known (see see \cite{BSh1}) than 
$$
{\rm C}_{{\rm A},\be}(\C_+)\otimes_{\rm h}{\rm C}_{{\rm A},\be}(\C_+)\subset\fM_{\rm A}(\C_+^2)
$$
and $\|\Phi\|_{\fM(\R^2)}\le\|\Phi\|_{{\rm C}_{{\rm A},\be}\otimes_{\rm h}{\rm C}_{{\rm A},\be}}$ for all
$\Phi\in{\rm C}_{{\rm A},\be}(\C_+)\otimes_{\rm h}{\rm C}_{{\rm A},\be}(\C_+)$,
see also \cite{PeHo} and \cite{APol}.

Let us now show that
$$
{\rm C}_{{\rm A},\be}(\C_+)\otimes_{\rm h}{\rm C}_{{\rm A},\be}(\C_+)
\subset{\rm C}^{\rm S}_{{\rm A},\be}(\C_+^2).
$$
Let $\Phi\in {\rm C}_{{\rm A},\be}(\C_+)\otimes_{\rm h}{\rm C}_{{\rm A},\be}(\C_+)$. Then $\Phi$
can be represented in the form $\Phi(z,w)=\sum_{n=1}^\be f_n(z)g_n(w)$ for $z,w\in\clos\C_+$,
where $f_n,g_n\in{\rm C}_{{\rm A},\be}(\C_+)$ for all $n\ge1$ such that
$\sum_{n=1}^\be |f_n|^2\le\const$ and $\sum_{n=1}^\be |g_n|^2\le\const$ everywhere on $\clos\C_+$.

Clearly, for each fixed $w\in\clos\C_+$, the series $\sum_{n=1}^\be f_n(\cdot)g_n(w)$ converges
uniformly on $\clos\C_+$. Hence, $\sum_{n=1}^\be f_n(\cdot)g_n(w)\in{\rm C}_{{\rm A},\be}(\C_+)$ for all $w\in\clos\C_+$

In the same way, $\sum_{n=1}^\be f_n(z)g_n(\cdot)\in{\rm C}_{{\rm A},\be}(\C_+)$ for all $z\in\clos\C_+$.

It remains  to prove that
$$
\fM_{\rm A}(\C_+^2)\cap{\rm C}_{{\rm A},\be}^{\rm S}(\C_+^2)\subset{\rm C}_{{\rm A},\be}(\C_+)\otimes_{\rm h}{\rm C}_{{\rm A},\be}(\C_+)
$$
and $\|\Phi\|_{{\rm C}_{{\rm A},\be}\otimes_{\rm h}{\rm C}_{{\rm A},\be}}\le\|\Phi\|_{\fM(\R^2)}$ for every
$\Phi\in{\rm C}_{{\rm A},\be}(\C_+)\otimes_{\rm h}{\rm C}_{{\rm A},\be}(\C_+)$.

Let $\Phi\in\fM_{\rm A}(\C_+^2)\cap{\rm C}_{{\rm A},\be}^{\rm S}(\C_+^2)$.
We have already observed that $\Phi\in\fM_{\rm d}(\R^2)$. By Theorem 5.1 of \cite{Pi}, $\Phi$ can be represented
in the form $\Phi(x,s)=(u_x,v_s)$,  where $\{u_x\}_{x\in\R}$  and $\{v_s\}_{s\in\R}$ 
are families in a Hilbert space $\h$ such that $\|u_x\|_\h^2\le\|\Phi\|_{\fM(\R^2)}$ and $\|v_s\|_\h^2\le\|\Phi\|_{\fM(\R^2)}$ for all $x$ and $s$ in $\R$. By the hypotheses of the theorem, $\Phi$ is continuous in each variable. Thus, by Theorem 2.2.3 of \cite{APol}, we can assume in addition that the maps $x\mapsto u_x$ and $s\mapsto v_s$
are weakly continuous and the linear span of each of the families $\{u_x\}_{x\in\R}$  and $\{v_s\}_{s\in\R}$
is dense in $\h$.
Hence, the space $\h$ is separable. Let $\{e_k\}_{k=1}^n$ be an orthonormal basis in $\h$,
where $n=\dim\h$. The vectors $u_x$ and $v_s$ can be represented in the form 
$u_x=\sum_{k=1}^n\f_k(x)e_k$ and $v_s=\sum_{k=1}^n\ov{\psi_k(s)}e_k$.

Clearly, $\f_k,\psi_k\in{\rm C}_{\rm b}(\R)$ for all $k$. Let us prove that $\f_k,\psi_k\in{\rm C}_{{\rm A},\be}$.

It is easy to see that the inclusion 
$\Phi\in{\rm C}^{\rm S}_{{\rm A},\be}(\C_+^2)$ implies that $\sum_{k=1}^n\f_k(x) \psi_k(s)\in{\rm C}^{\rm S}_{{\rm A},\be}(\C_+^2)$.
Then $\sum_{k=1}^n\f_k(\cdot) \psi_k(s)\in{\rm C}_{{\rm A},\be}$ for all $s\in\R$. Hence,
$\sum_{k=1}^n\a_k\f_k\in{\rm C}_{{\rm A},\be}$ for every sequence $\{\a_k\}_{k=1}^n$ of complex numbers
such that $\sum_{k=1}^n|\a_k|^2<\be$. In particular, $\f_k\in{\rm C}_{{\rm A},\be}$ for all $k$.
In the same way $\psi_k\in{\rm C}_{{\rm A},\be}$ for all $k$.
It remains to prove that $\Phi(z,w)=\sum_{k=1}^n\f_k(z)\psi_k(w)$ for all $z,w\in\clos \C_+$. 
In particular, this equality holds for $z,w\in\R$. It remains to observe that $\Phi$ is a bounded, continuous in each
variable and analytic in each variable function on $\C_+$. 
$\bl$

\medskip

We proceed now to the description of analytic Schur multipliers in ${\rm C}^{\rm S}_{\rm A}(\C_+^2)$.

\begin{thm}
\label{kon}
$$
\fM_{\rm A}(\C_+^2)\cap{\rm C}_{\rm A}^{\rm S}(\C_+^2)=\CAr\otimes_{\rm h}\CAr
$$
and  
$$
\|\Phi\|_{\CAr\otimes_{\rm h}\CAr}=\|\Phi\|_{\fM(\R^2)}\quad
\mbox{for every}\quad\Phi\in\CAr\otimes_{\rm h}\CAr.
$$
\end{thm}

\Pf In the same way as in the proof of Theorem \ref{beskon} we obtain
$$
\CAr\otimes_{\rm h}\CAr\subset\fM_{\rm A}(\C_+^2)\cap{\rm C}_{\rm A}^{\rm S}(\C_+^2)
$$
and $\|\Phi\|_{\fM(\R^2)}\le\|\Phi\|_{\CAr\otimes_{\rm h}\CAr}$

It remains to prove that
$$
\fM_{\rm A}(\C_+^2)\cap{\rm C}_{\rm A}^{\rm S}(\C_+^2)\subset\CAr\otimes_{\rm h}\CAr
$$
and $\|\Phi\|_{\CAr\otimes_{\rm h}\CAr}\le\|\Phi\|_{\fM(\R^2)}$.

Repeating the same arguments as in the proof of Theorem \ref{beskon},
we can represent $\Phi \big| \R^2$ in the form $\Phi(x,s)=(u_x,v_s)$, where $\{u_x\}_{x\in\R}$  and $\{v_s\}_{s\in\R}$
are families in a Hilbert space $\h$ such that the maps $x\mapsto u_x$ and $s\mapsto v_s$
are weakly continuous, $\|u_x\|_\h^2\le\|\Phi\|_{\fM(\R^2)}$ and $\|v_s\|_\h^2\le\|\Phi\|_{\fM(\R^2)}$
for all $x,s\in\R$. Moreover,  the linear span of each of the families $\{u_x\}_{x\in\R}$  and $\{v_s\}_{s\in\R}$
is dense in $\h$. Thus, the space $\h$ is separable. Let $\{e_k\}_{k=1}^n$ be an orthonormal basis in $\h$,
where $n=\dim\h$. The vectors $u_x$ and $v_s$ can be represented in the form
$u_x=\sum_{k=1}^n\f_k(x)e_k$ and $v_s=\sum_{k=1}^n\ov{\psi_k(s)}e_k$.

Clearly, $\f_k,\psi_k\in{\rm C}(\widehat\R)$ for all $k$. Let us prove that
$\f_k$ and $\psi_k$ belong to $\CAr$.

It is easy to see that the inclusion
$\Phi\in{\rm C}_{\rm A}(\C_+^2)$ implies that $\sum_{k=1}^n\f_k(x) \psi_k(s)\in{\rm C}^{\rm S}_{\rm A}(\C_+^2)$.
Then $\sum_{k=1}^n\f_k(\cdot) \psi_k(s)\in{\rm C}_{\rm A}(\C_+)$ for all $s\in\R$. Hence,
$\sum_{k=1}^n\a_k\f_k\in{\rm C}_{\rm A}(\C_+)$ for every sequence $\{\a_k\}_{k=1}^n$ of complex numbers
such that $\sum_{k=1}^n|\a_k|^2<\be$. In particular, $\f_k\in{\rm C}_{\rm A}(\C_+)$ for all $k$.
In the same way, $\psi_k\in{\rm C}_{\rm A}(\C_+)$ for all $k$.
We have to prove that $\Phi(z,w)=\sum_{k=1}^n\f_k(z)\psi_k(w)$ for all $z,w\in\clos \C_+$.
This equality holds for $z,w\in\R$. It remains to observe that $\Phi$ is a bounded, continuous in each
variable and analytic in each variable function on $\C_+$.
$\bl$

\

\section{\bf A representation of the class of analytic Schur multipliers as a Haagerup tensor product}
\label{vsyudu}
\setcounter{equation}{0}

\

In this section we describe the space of analytic Schur multipliers $\fM_{\rm A}(\C_+^2)$ as the
Haagerup tensor product $H^\be_{\rm r}(\C_+)\otimes_{\rm h}{H^\be_{\rm r}}(\C_+)$.
       
Let $\Phi\in H^\be(\C_+^2)$. Put $\Phi_t(z,w)=\Phi(z+t\ri,w+t\ri)$, where $t>0$ and $z,w\in\clos\C_+$.

\begin{thm}
\label{dandh}
$$
\fM_{\rm A}(\C_+^2)=H^\be_{\rm r}(\C_+)\otimes_{\rm h}{H^\be_{\rm r}}(\C_+)
$$
and $\|\Phi\|_{H^\be_{\rm r}(\C_+)\otimes_{\rm h}{H^\be_{\rm r}}(\C_+)}=\|\Phi\|_{\fM_{\rm d}(\R^2)}$
for all $\Phi \in\fM_{\rm A}(\C_+^2)$.
\end{thm}

\medskip

We need the following lemma.

\begin{lem}
\label{AM}
Let $\Phi \in \fM_{\rm A}(\C_+^2)$.
Then $\Phi\in\fM_{\rm d}(\C_+^2)$ and $\|\Phi\|_{\fM_{\rm d}(\C_+^2)}\le\|\Phi\|_{\fM_{\rm d}(\R^2)}$.
\end{lem}

\Pf If we assume in addition that $\Phi\in{\rm C}^{\rm S}_{{\rm A},\be}(\C_+^2)$, the result follows from Theorem \ref{beskon}.

Let $\Phi\in  H^\be_{\rm r}(\C_+^2)$. Consider the Poisson kernel $P_t$ defined by $P_t(x)\df\frac1\pi\cdot\frac t{x^2+t^2}$, $x\in\R$.

Put $\bs{P}_{t}(x,y)\df P_t(x)P_t(y)$. Then $\Phi_{t}=\Phi*\bs{P}_t$. Since the norm in $\fM_{\rm d}(\R^2)$ is translation invariant, it follows that
$$
\|\Phi_{t}\|_{\fM_{\rm d}(\R^2)}\le\|\Phi\|_{\fM_{\rm d}(\R^2)}\|\bs{P}_t\|_{L^1(\R^2)}=\|\Phi\|_{\fM_{\rm d}(\R^2)}.
$$
Applying Theorem \ref{beskon}, we find that $\|\Phi_{t}\|_{\fM_{\rm d}(\clos\C_+^2)}\le\|\Phi\|_{\fM_{\rm d}(\R^2)}$
for all $t>0$.
It remains to observe that $\|\Phi\|_{\fM_{\rm d}(\C_+^2)}=\lim_{t\to0^+}\|\Phi_t\|_{\fM_{\rm d}(\clos\C_+^2)}$. $\bl$

\medskip

 {\bf Proof of Theorem \ref{dandh}} In the same way as in the proofs of Theorem \ref{beskon} and Theorem \ref{kon} we obtain
$$
H^\be_{\rm r}(\C_+)\otimes_{\rm h}{H^\be_{\rm r}}(\C_+)\subset\fM_{\rm A}(\C_+^2)
$$
and $\|\Phi\|_{\fM_{\rm d}(\R^2)}\le\|\Phi\|_{H^\be_{\rm r}(\C_+)\otimes_{\rm h}{H^\be_{\rm r}}(\C_+)}$.

It remains to prove that
$$
\fM_{\rm A}(\C_+^2)\subset H^\be_{\rm r}(\C_+)\otimes_{\rm h}{H^\be_{\rm r}}(\C_+)
$$
and $\|\Phi\|_{H^\be_{\rm r}(\C_+)\otimes_{\rm h}{H^\be_{\rm r}}(\C_+)}\le\|\Phi\|_{\fM_{\rm d}(\R^2)}$.

Let $\Phi\in\fM_{\rm A}(\C_+^2)$.
By Lemma \ref{AM}, we have
 $\|\Phi\|_{\fM_{\rm d}(\C_+^2)}\le\|\Phi\|_{\fM_{\rm d}(\R^2)}$. Besides, 
$$\|\Phi\|_{\fM_{\rm d}(\R^2)}\le\|\Phi\|_{\fM_{\rm d}(\clos\C_+^2)}\le\|\Phi\|_{H^\be_{\rm r}(\C_+)\otimes_{\rm h}{H^\be_{\rm r}}(\C_+)}.
$$
 
It remains to prove that $\Phi\in H^\be_{\rm r}(\C_+)\otimes_{\rm h}{H^\be_{\rm r}}(\C_+)$ and
$\|\Phi\|_{H^\be_{\rm r}(\C_+)\otimes_{\rm h}{H^\be_{\rm r}}(\C_+)}\le\|\Phi\|_{\fM_{\rm d}(\R^2)}$.

By Theorem 5.1 of \cite{Pi}, $\Phi|\C_+^2$ can be represented
in the form $\Phi(z,w)=(u_z,v_w)$,  where $\{u_z\}_{z\in\C_+}$  and $\{v_w\}_{z\in\C_+}$
are families in a Hilbert space $\h$ such that $\|u_z\|_\h^2\le\|\Phi\|_{\fM_{\rm d}(\R^2)}$ and
$\|v_w\|_\h^2\le\|\Phi\|_{\fM_{\rm d}(\R^2)}$ for all $z$ and $w$ in $\C_+$. By the hypotheses of the theorem, $\Phi|\C_+^2$ is continuous in each variable. Thus, by Theorem 2.2.3 of \cite{APol}, we can assume in addition that the maps $z\mapsto u_z$ and $w\mapsto v_w$
are weakly continuous on $\C_+$ and the linear spans of both families $\{u_z\}_{z\in\C_+}$  and $\{v_w\}_{w\in\C_+}$ are dense in $\h$.

It follows from the facts that the functions $z\mapsto u_z$ and $w\mapsto v_w$ are weakly continuous in
$\C_+$ that the space $\h$ is separable.

Let $z\in\R$ and $w\in\C_+$.  Then we have $\Phi(z,w)=\lim_{t\to 0^+}\Phi(z+t\ri,w)$
for all $w\in\C_+$, i.e. the limit $\lim_{t\to 0^+}(u_{z+t\ri},v_w)$ exists for all $w\in\C_+$.
Put $u_z\df\lim_{t\to 0^+}u_{z+t\ri}$ (the limit exists in the weak topology of $\h$).
In the same way for $w\in\R$, we can define the vector $v_w\df\lim_{t\to 0^+}v_{w+t\ri}$ (the limit exists in the weak topology of $\h$).
Then  $\Phi(z,w)=(u_z,v_w)$ for all $z\in\C_+$
Thus, we have $\Phi(z,w)=(u_z,v_w)$ for all $z,w\in\clos\C_+$ if $z\notin\R$ or $w\notin\R$.
If $z,w\in\R$, then we have $\Phi(z,w)=\lim_{t\to 0^+}(u_{z+t\ri},v_w)=(u_z,v_w)$
because $u_z=\lim_{t\to 0^+}u_{z+t\ri}$ in the weak topology of $\h$.

Let $\{e_k\}_{k=1}^n$ be an orthonormal basis in $\h$,
where $n=\dim\h$. The vectors $u_z$ and $v_w$ can represented in the form
$u_z=\sum_{k=1}^n\f_k(z)e_k$ and $v_w=\sum_{k=1}^n\ov{\psi_k(w)}e_k$,
where $\f_k$ and $\psi_k$ are bounded analytic functions on $\C_+$.
It remains to prove that $\f_k, \psi_k\in H^\be_{\rm r}(\C_+)$ for all $k$.

Put $\h_0=\{h\in\h: (u_z,h)\in H^\be_{\rm r}(\C_+)\}$ and $\h_1=\{h\in\h: (h,v_w)\in H^\be_{\rm r}(\C_+)\}$.
Clearly, $\h_0$ and $\h_1$ are closed subspaces $\h$. Moreover, $v_w\in\h_0$ for all $w\in\clos\C_+$
and $u_z\in\h_1$ for all $z\in\clos\C_+$. Hence, $\h_0=\h_1=\h$. It remains to observe
that $\f_k(z)=(u_z,e_k)$ and $\psi_k(w)=(e_k,v_w)$. $\bl$

\

\pagebreak

\section{\bf Analytic Schur multipliers in the case of the bidisc}
\label{bid}
\setcounter{equation}{0}

\

In this section we are going to show that using the same techniques,
one can obtain similar results for
functions analytic in the bidisc.

We need the following analog of the space $H^\be_{\rm r}(\C_+^2)$:

We denote by $H^\be_{\rm r}(\dd^2)$ the space of functions $\Phi$ on $\clos\dd^2$
such that $\Phi|\dd^2\in H^\be(\dd^2)$ and $\lim_{r\to 1^-}\Phi(rz,w)=\lim_{r\to 1^-}\Phi(z,rw)=\Phi(z,w)$
for arbitrary $z,w\in\clos\dd$.

Put
$$
{\rm C}_{\rm A}^{\rm S}(\dd^2)=\big\{f\in H^\be(\dd^2):~f(z,\cdot), ~\:f(\cdot,z)\in{\rm C}_{\rm A}\quad\mbox{for all}\quad z\in\clos\dd\big\},
$$
where ${\rm C}_{\rm A}$ stands for the disc algebra.

As in the case of the classes of functions on $\C_+^2$, the restrictions of functions  in 
${\rm C}^{\rm S}_{{\rm A}}(\dd_+^2)$ to $\T^2$ do not have to be continuous functions on $\T^2$, while
the restrictions of functions in ${\rm C}^{\rm S}_{{\rm A}}(\dd_+^2)$  to $\dd_+^2$ must be holomorphic in 
$\dd_+^2$.

Let $\fM_{\rm A}(\dd^2)$ be the class of {\it analytic Schur multipliers} on $\dd^2$, i.e. the class of
functions  $f\in H^\be_{\rm r}(\dd^2)$ such that $f|\T^2$ belongs to the space $\fM(\T^2)$ of {\it Borel Schur multipliers}, i.e. the space of functions on $\T\times\T$ that are Schur multipliers with respect to arbitrary Borel spectral measures on $\T$.

The following theorem can be proved in the same way as Theorem \ref{kon}.
Moreover, it can be deduced from Theorem \ref{kon} by applying the conformal map of the unit disc onto the upper half-plane.

\begin{thm}
\label{beskondisk}
$$
\fM_{\rm A}(\dd^2)\cap{\rm C}_{\rm A}^{\rm S}(\dd^2)=\rm C_{\rm A}\otimes_{\rm h}{\rm C}_{\rm A}
$$
Let $\Phi\in{\rm C}_{\rm A}^{\rm S}(\dd^2)\cap\fM_{\rm A}(\dd^2)$.
Then $\Phi\in{\rm C}_{\rm A}\otimes_{\rm h}{\rm C}_{\rm A}$ and
$$
\|\Phi\|_{\rm C_{\rm A}\otimes_{\rm h}{\rm C}_{\rm A}}
=\|\Phi\|_{\fM(\T^2)}.
$$
\end{thm}

\medskip

\medskip

Note that in the special case when $\Phi$ is the divided difference of a function on $\T$, the conclusion
of Theorem \ref{beskondisk} was obtained in \cite{KS}. Let us also mention that
in \cite{A} the result of \cite{KS} was generalized to the case when
$\Phi$ is the divided difference of a function on a closed subset of $\C$ without isolated points.

Now we proceed to an analog of Theorem \ref{dandh} for the disc. We start with a lemma.

Let $\Phi\in H^\be(\dd^2)$. Put $\Phi_r(z,w)=\Phi(rz,rw)$, where $r\in[0,1)$ and $z,w\in\clos\dd_+$.

\begin{lem}
\label{AMdisk}
Let $\Phi \in \fM_{\rm A}(\dd^2)$.
Then $\Phi\in\fM_{\rm d}(\dd^2)$ and $\|\Phi\|_{\fM_{\rm d}(\dd^2)}\le\|\Phi\|_{\fM_{\rm d}(\T^2)}$.
\end{lem}

\Pf If we assume in addition that $\Phi\in{\rm C}_{\rm A}^{\rm S}(\dd^2)$, the result follows from Theorem \ref{beskondisk}.

Suppose now that $\Phi$ is an arbitrary function in $H^\be_{\rm r}(\dd^2)$. Consider the Poisson kernel
${\mathcal P}_r$ defined by ${\mathcal P}_r(\z)\df\frac{1-r^2}{|r-\z|^2}$, $\z\in\T$,
$r\in[0,1)$.

Put $\bs{{\mathcal P}}_r(\z,\xi)\df {\mathcal P}_r(\z){\mathcal P}_r(\xi)$. Then $\Phi_r=\Phi*\bs{{\mathcal P}}_r$. Since the norm in
$\fM_{\rm d}(\T^2)$ is rotation invariant, it follows that
$$
\|\Phi_r\|_{\fM_{\rm d}(\T^2)}\le\|\Phi\|_{\fM_{\rm d}(\T^2)}\|{\mathcal P}_r\|_{L^1(\T^2,\m_2)}=\|\Phi\|_{\fM_{\rm d}(\T^2)}.
$$
Applying Theorem \ref{beskondisk}, we find that $\|\Phi_r\|_{\fM_{\rm d}(\clos\dd^2)}\le\|\Phi\|_{\fM_{\rm d}(\T^2)}$
for all $r\in[0,1)$.
It remains to observe that $\|\Phi\|_{\fM_{\rm d}(\clos\dd^2)}=\lim_{r\to1^-}\|\Phi_r\|_{\fM_{\rm d}(\clos\dd^2)}$. $\bl$

%
%
%

\begin{thm}
\label{nikakaya}
$$
\fM_{\rm A}(\dd^2)=H^\be_{\rm r}(\dd)\otimes_{\rm h}{H^\be_{\rm r}}(\dd)
$$
and $\|\Phi\|_{H^\be_{\rm r}(\dd)\otimes_{\rm h}{H^\be_{\rm r}}(\dd)}=\|\Phi\|_{\fM_{\rm A}(\dd^2)}$
for all $\Phi \in \fM_{\rm A}(\dd^2)$.
\end{thm}

\Pf  In the same way as in the proofs of Theorem \ref{dandh} we have
$$
H^\be_{\rm r}(\dd)\otimes_{\rm h}{H^\be_{\rm r}}(\dd)\subset H^\be_{\rm r}(\dd^2)
$$
and  $\|\Phi\|_{\fM_{\rm d}(\T^2)}\le\|\Phi\|_{H^\be_{\rm r}(\dd)\otimes_{\rm h}{H^\be_{\rm r}}(\dd)}$.

It remains to show that
$$
H^\be_{\rm r}(\dd^2)\subset H^\be_{\rm r}(\dd)\otimes_{\rm h}{H^\be_{\rm r}}(\dd)
$$
and $\|\Phi\|_{H^\be_{\rm r}(\dd)\otimes_{\rm h}{H^\be_{\rm r}}(\dd)}\le\|\Phi\|_{\fM_{\rm d}(\T^2)}$.

Let $\Phi\in\fM_{\rm A}(\dd^2)$.
By Lemma \ref{AMdisk}, we have  $\|\Phi\|_{\fM_{\rm d}(\clos\dd^2)}=\|\Phi\|_{\fM_{\rm d}(\T^2)}$.
It suffices to show that $\Phi\in H^\be_{\rm r}(\dd)\otimes_{\rm h}{H^\be_{\rm r}}(\dd)$ and
$\|\Phi\|_{H^\be_{\rm r}(\dd)\otimes_{\rm h}{H^\be_{\rm r}}(\dd)}\le\|\Phi\|_{\fM_{\rm d}(\T^2)}$.

Thus, by Theorem 5.1 of \cite{Pi}, $\Phi$ can be represented
in the form $\Phi(z,w)=(u_z,v_w)$,  where $\{u_z\}_{z\in\dd}$  and $\{v_w\}_{w\in\dd}$
are families in a Hilbert space $\h$ such that $\|u_z\|_\h^2\le\|\Phi\|_{\fM_{\rm d}(\T^2)}$ and
$\|v_w\|_\h^2\le\|\Phi\|_{\fM_{\rm d}(\T^2)}$ for all $z$ and $w$ in $\dd$. By the hypotheses of the theorem, $\Phi|\T^2$ is continuous in each variable. Thus, by Theorem 2.2.3 of \cite{APol}, we can assume in addition that the maps $z\mapsto u_z$ and $w\mapsto v_w$
are weakly continuous on $\dd$ and the linear span of each of the families $\{u_z\}_{z\in\dd}$  and $\{v_w\}_{w\in\dd}$
is dense in $\h$.

It follows from the facts that the functions $z\mapsto u_z$ and $w\mapsto v_w$ are weakly continuous in
$\dd$ that the space $\h$ is separable.

Let $z\in\T$ and $w\in\dd$.  Then we have $\Phi(z,w)=\lim_{r\to 1^-}\Phi(rz,w)$
for all $w\in\dd$, i.e. the limit $\lim_{r\to 1^-}(u_{rz},v_w)$ exists for all $w\in\dd$.
Put $u_z\df\lim_{r\to 1^-}u_{rz}$ (the limit exists in the weak topology of $\h$).
In the same way for $w\in\T$, we can define the vector $v_w\df\lim_{r\to 1^-}v_{rw}$ (the limit exists in the weak topology of $\h$).
Then  $\Phi(z,w)=(u_z,v_w)$ for all $z\in\C_+$
Thus, we have $\Phi(z,w)=(u_z,v_w)$ for all $z,w\in\clos\dd$ if $z\notin\T$ or $w\notin\T$.
If both $z$ and $w$ belong to $\T$, then we have $\Phi(z,w)=\lim_{r\to 1^-}(u_{rz},v_w)=(u_z,v_w)$
because $u_z=\lim_{r\to 1^-}u_{rz}$ in the weak topology of $\h$.

Let $\{e_k\}_{k=1}^n$ be an orthonormal basis in $\h$,
where $n=\dim\h$. The vectors $u_z$ and $v_w$ can be represented in the form
$u_z=\sum_{k=1}^n\f_k(z)e_k$ and $v_w=\sum_{k=1}^n\ov{\psi_k(w)}e_k$,
where $\f_k$ and $\psi_k$ are bounded analytic functions in $\dd$.
It remains to prove that $\f_k, \psi_k\in H^\be_{\rm r}(\dd)$ for all $k$.

Put $\h_0=\{h\in\h: (u_z,h)\in H^\be_{\rm r}(\dd)\}$ and $\h_1=\{h\in\h: (h,v_w)\in H^\be_{\rm r}(\dd)\}$.
Clearly, $\h_0$ and $\h_1$ are closed subspaces $\h$. Moreover, $v_w\in\h_0$ for all $w\in\clos\dd$
and $u_z\in\h_1$ for all $z\in\clos\dd$. Hence, $\h_0=\h_1=\h$. It remains to observe
that $\f_k(z)=(u_z,e_k)$ and $\psi_k(w)=(e_k,v_w)$. $\bl$

\

\begin{footnotesize}

\
 
\noindent
\begin{tabular}{p{8cm}p{15cm}}
A.B. Aleksandrov & V.V. Peller \\
St.Petersburg State University & St.Petersburg State University \\
Universitetskaya nab., 7/9  & Universitetskaya nab., 7/9\\
199034 St.Petersburg, Russia & 199034 St.Petersburg, Russia \\
\\

St.Petersburg Department &St.Petersburg Department\\
Steklov Institute of Mathematics  &Steklov Institute of Mathematics  \\
Russian Academy of Sciences  & Russian Academy of Sciences \\
Fontanka 27, 191023 St.Petersburg &Fontanka 27, 191023 St.Petersburg\\
Russia&Russia\\
email: alex@pdmi.ras.ru&email: pellerv@gmail.com

\end{tabular}
\end{footnotesize}

\end{document}